\documentclass[12pt]{amsart}
\usepackage[]{hyperref} 
\usepackage{amssymb}
\usepackage[curve]{xypic}
\newcommand{\Omit}[1]{\begin{tiny}#1\end{tiny}}
\renewcommand{\Omit}[1]{}

\newbox\mybox
\def\overtag#1#2#3{\setbox\mybox\hbox{$#1$}\hbox to
  0pt{\vbox to 0pt{\vglue-#3\vglue-\ht\mybox\hbox to \wd\mybox
      {\hss$\scriptstyle#2$\hss}\vss}\hss}\box\mybox}
\def\undertag#1#2#3{\setbox\mybox\hbox{$#1$}\hbox to 0pt{\vbox to
    0pt{\vglue#3\vglue\ht\mybox\hbox to \wd\mybox
      {\hss$\scriptstyle#2$\hss}\vss}\hss}\box\mybox}
\def\lefttag#1#2#3{\hbox to 0pt{\vbox to 0pt{\vss\hbox to
      0pt{\hss$\scriptstyle#2$\hskip#3}\vss}}#1}
\def\righttag#1#2#3{\hbox to 0pt{\vbox to 0pt{\vss\hbox to
      0pt{\hskip#3$\scriptstyle#2$\hss}\vss}}#1}

\def\Dot{\lower.2pc\hbox to 2.5pt{\hss$\bullet$\hss}}
\def\Circ{\lower.2pc\hbox to 2.5pt{\hss$\circ$\hss}}
\def\Vdots{\raise5pt\hbox{$\vdots$}}
\def\splicediag#1#2{\xymatrix@R=#1pt@C=#2pt@M=0pt@W=0pt@H=0pt}

\renewcommand\frame[2][3pt]{\hbox{$\vcenter{\hbox{\vrule\vbox
{\hrule\kern#1\hbox{\kern#1$#2$\kern#1}\kern#1\hrule}\vrule}}$}}

\newcommand\lineto{\ar@{-}}
\newcommand\dashto{\ar@{--}}
\newcommand\dotto{\ar@{.}}
\newcommand{\bt}{\bullet}
\newcommand{\C}{\mathbb C}

\newcommand{\Z}{\mathbb Z}
\newcommand{\Q}{\mathbb Q}

\newtheorem*{ECtheorem}{Main Conjecture}
\newtheorem*{RatConj}{Rational Conjecture}
\newtheorem*{QGor}{$\Q$-Gorenstein Conjecture}
\newtheorem*{Codim}{Codimension 2 Conjecture}

\newtheorem*{theorem*}{Theorem}
\newtheorem*{corollary*}{Corollary}
\newtheorem{theorem}{Theorem}[section]
\newtheorem{proposition}[theorem]{Proposition}
\newtheorem{lemma}[theorem]{Lemma}
\newtheorem{corollary}[theorem]{Corollary}
\newtheorem{conjecture}[theorem]{Conjecture}
\theoremstyle{definition}
\newtheorem{example}[theorem]{Example}

\newtheorem*{example*}{Example}

\newtheorem{remark}[theorem]{Remark}
\newtheorem*{remark*}{Remark}

\newtheorem*{conjecture*}{Conjecture}
\begin{document}
\title{Milnor and Tjurina numbers for smoothings of surface singularities}
\author{Jonathan Wahl}
\dedicatory { To Eduard Looijenga on his 65th birthday}
\address{Department of Mathematics\\The University of North
  Carolina\\Chapel Hill, NC 27599-3250} \email{jmwahl@email.unc.edu}
  \begin{abstract} For an isolated hypersurface singularity $\{f=0\}$, the Milnor number $\mu$ is greater than or equal to the Tjurina number $\tau$ (the dimension of the base of the semi-universal deformation), with equality if $f$ is quasi-homogeneous.  K. Saito proved the converse.  The same result is true for complete intersections, but is much harder.  For a Gorenstein surface singularity $(V,0)$,  the difference $\mu-\tau$ can be defined whether or not  $(V,0)$ is smoothable; the author has proved it is non-negative, and equal to $0$ if and only if $(V,0)$ is quasi-homogeneous.  We conjecture a similar result for non-Gorenstein surface singularities.  Here,  $\mu-\tau$ must be modified so that it is independent of any smoothing.  This expression, involving cohomology of exterior powers of the bundle of logarithmic derivations on the minimal good resolution, is conjecturally non-negative, with equality iff $(V,0)$ is quasi-homogeneous.  We prove the ``if" part; identify special cases where the conjecture is particularly interesting; verify it in some non-trivial cases; and prove it for a $\Q$-Gorenstein smoothing when the index one cover is a hypersurface.  This conjecture arose regarding the classification of surface singularities with rational homology disk smoothings.
    \end{abstract}
\maketitle

\section{Introduction}

Suppose the $n$-dimensional hypersurface $V=\{f(z_1,...,z_{n+1})=0\}\subset \C^{n+1}$ has an isolated singularity at the origin.  Then the Milnor fibre $M$ of $(V,0)$ has the homotopy type of a bouquet of a certain number $\mu$ of $n$-spheres, where the Milnor number $\mu$ is the length of the Jacobian algebra $$J_f=\C[[z_1,..., z_{n+1}]]/(\partial f/ \partial z_1,...,\partial f/\partial z_{n+1}).$$  The base space of the semi-universal deformation of $(V,0)$ has dimension $\tau,$ which is the length of the Tjurina algebra $$T_f=\C[[z_1,..., z_{n+1}]]/(f,\partial f/ \partial z_1,...,\partial f/\partial z_{n+1}).$$
Clearly, $\mu \geq \tau,$ with equality iff $f$ belongs to the Jacobian ideal.  This is the case  when $f$ is  a weighted homogeneous polynomial.  K. Saito \cite{saito} proved the converse:  equality implies that $f$ is quasi-homogeneous, i.e., analytically equivalent to such a polynomial.

A local complete intersection singularity $(V,0)$ of dimension $n\geq 1$ also has a Milnor fibre $M$ which is a bouquet of a certain number $\mu$ of $n$-spheres.  The base space of its semi-universal deformation is again smooth of dimension $\tau=l(\text{Ext}^1_R(\Omega^1_R,R))$, where $l$ denotes length and $R$ is the local ring $\mathcal O_{V,0}$.  In this case it is harder to relate the size of the two invariants, which are dimensions of very different-looking spaces.  But it was eventually proved by work of (among others) Greuel \cite{greuel}, Looijenga-Steenbrink \cite{LS}, the author \cite{w}, and Vosegaard \cite{vose} that also in these cases one has $\mu \geq \tau$, with equality if and only if $(V,0)$ is quasi-homogeneous.

Smoothings for $(V,0)$ a general normal surface singularity were investigated in \cite{wahls}.  A singularity (e.g., the cone over the rational quartic curve in $\mathbb P^4$) can have several topologically distinct smoothings,  occurring over ``smoothing components" of different dimension in the base space of the semi-universal deformation of $(V,0)$.  Thus, for a smoothing $\pi:(\mathcal V,0)\rightarrow (\mathbb C,0)$, a Milnor fibre $M$ may be defined, but the rank $\mu_{\pi}$ of $H^2(M)$, and the dimension $\tau_{\pi}$ of the corresponding smoothing component, depend on the smoothing.  (It is known that the first Betti number of the Milnor fibre is $0$, by \cite{gst}).  Define  $\alpha_{\pi}$ to be the colength of the restriction map of the dual of (relative) dualizing sheaves $$\omega_{\mathcal V/\C}^*\otimes \mathcal O_V\rightarrow \omega_V^*.$$   The following theorem was proved in special cases in \cite{wahls}, and in general modulo several conjectures which were later established by Greuel, Looijenga, and Steenbrink (\cite{gst}, \cite{gl}, \cite{lo}).

\begin{theorem}[\cite{wahls}] Let $\pi:(\mathcal V,0)\rightarrow (\C,0)$ be a smoothing of a normal surface singularity $(V,0)$, with $\mu_{\pi}, \tau_{\pi}, \alpha_{\pi}$ as above.  Let $(X,E)\rightarrow (V,0)$ be a good resolution.
Then 
\begin{enumerate}
\item $1+\mu_{\pi}=\alpha_{\pi} +13h^1(\mathcal O_X)+\chi_T(E)-h^1(-K_X).$
\medskip

\item \ \ \ $\tau_{\pi}=2\alpha_{\pi} +12h^1(\mathcal O_X)+h^1(\Theta_X)-2h^1(-K_X).$
\end{enumerate}
 If $(V,0)$ is Gorenstein, then $\alpha_{\pi}=0$, so $\mu$ and $\tau$ are independent of the smoothing.
\end{theorem}

In the Gorenstein case, one may use the expressions from the preceding theorem to define (possibly negative) singularity invariants $\mu$ and $\tau$, independent of smoothability.  With this definition, we have in general

\begin{theorem}[\cite{w}] If $(V,0)$ is a Gorenstein surface singularity, then $\mu- \tau \geq 0$, with equality if and only if $(V,0)$ is weighted homogeneous.
\end{theorem}


The point of this work is to study non-Gorenstein (e.g., rational) surface singularities.   Rather than $\mu-\tau$, consider the intrinsic invariant $$(\mu_{\pi}-\tau_{\pi})+\alpha_{\pi},$$
depending only on $(V,0)$.  A more useful version of Theorem 1.1 is given in  Section 2:

\begin{corollary*}[\ref{s}]  Let $(X,E)\rightarrow (V,0)$ be the minimal good resolution of a normal surface singularity, not a rational double point.  Denote by $S_X=(\Omega^1_X(\text{log}(E))^*$ the sheaf of derivations on $X$, logarithmic along $E$.  Then for any smoothing $\pi$, 
$$1+(\mu_{\pi}-\tau_{\pi})+\alpha_{\pi}=h^1(\mathcal O_X)-h^1(S_X)+h^1(\wedge^2 S_X).$$
\end{corollary*}

In the Gorenstein case, $\alpha_{\pi}=0$ and the methods of \cite{w} allow one to compute the right hand side, particularly the difficult term $h^1(S_X)$ (whose dimension can vary in an ``equisingular family").  Note $h^1(\wedge^2 S_X)=h^1(-(K_X+E))$ is the \emph{second plurigenus} $\delta_2(V)$ of K. Watanabe (\cite{wa}, see (\ref{ok}) below), and it is not so easy to compute even for rational or quasi-homogeneous singularities.

The main purpose of this paper is to offer the following

\begin{ECtheorem}  Let $(X,E)\rightarrow (V,0)$ be the minimal good resolution of a non-Gorenstein normal surface singularity.  Denoting by $S_X$ the sheaf of logarithmic derivations on $X$, one has
$$h^1(\mathcal O_X)-h^1(S_X)+h^1(\wedge^2 S_X)\geq 0,$$ with equality if and only if $(V,0)$ is quasi-homogeneous.

\end{ECtheorem}
 The relevant expression $\sum_{i=0}^2 (-1)^i h^1(\wedge^i S_X)$ looks somewhat like a second Chern class, as explained in (\ref{ch}).
 
The Main Conjecture (given also in (\ref{conj}) below) might be overly optimistic.  One can check it in certain cases by computing the right or left side in Corollary \ref{s}.  

One implication is proved in this paper:

\begin{theorem*}[\ref{d}] If $(V,0)$ is quasi-homogeneous and not Gorenstein, then
$$h^1(\mathcal O_X)-h^1(S_X)+h^1(\wedge^2 S_X)= 0.$$ 
\end{theorem*}

Certain special cases are worth pursuing.  In Section $4$, we have the

\begin{RatConj}   Let $(X,E)\rightarrow (V,0)$ be the minimal good resolution of a rational surface singularity, not an RDP (rational double point).  Then $$h^1(S_X) \leq h^1(-(K_X+E)),$$ with equality if and only if $(V,0)$ is quasi-homogeneous.
\end{RatConj} 

Recall that $h^1(S_X)$ is the dimension of the smooth space of equisingular deformations of $(V,0)$, obtained from the deformations of $X$ which preserve every exceptional curve \cite{wan}.  As for $h^1(-(K_X+E))$, we prove in Corollary \ref{h1} that it can be computed from the resolution graph; thus the Conjecture would give a topological upper bound for the dimension $h^1(\Theta_X)$ of the Artin component of $(V,0)$.  As evidence for the Rational Conjecture, we show in Section 4 validity in several cases:
\begin{enumerate}
\item  $h^1(S_X)=0$ (\ref{a}).
\item the resolution graph of $(V,0)$ is star-shaped (\ref{b}).
\item an example  with non-star-shaped graph and $h^1(S_X)=1$ (\ref{c}).
\end{enumerate}
In a forthcoming paper, we verify the Conjecture (in a stronger form) for any rational graph which is ``sufficiently negative at the nodes''; see (\ref{jon}) for a precise statement.

Recall that $(V,0)$ is called $\Q$-Gorenstein if $K_V$ is $\Q$-Cartier, i.e. some $rK_V$ is invertible; then $(V,0)$ is an $r$-cyclic quotient of its index one (or canonical) cover $(W,0)$, which is Gorenstein.  A smoothing $\pi:(\mathcal V,0)\rightarrow (\C,0)$ is called $\Q$-Gorenstein if it is an $r$-cyclic quotient of a smoothing of $(W,0)$.   In this case, $\alpha_{\pi}=0$ (see Lemma \ref{alpha}), so the Conjecture is again about $\mu-\tau$.  We discuss in Section $5$ the

\begin{QGor} For a $\Q$-Gorenstein smoothing of a non-Gorenstein singularity $(V,0)$, $$\mu\geq \tau -1,$$ with equality if and only if $(V,0)$ is quasi-homogeneous.
\end{QGor}

We verify this last Conjecture in the case that the canonical cover is a hypersurface singularity.  In fact, a much more general result, in all dimensions, is proved.

\begin{theorem*}[\ref{q}] Let $(W,0)$ be an isolated hypersurface singularity $(\{f(z_1,\cdots,z_{n+1})=0\},0)\subset(\C^{n+1},0)$, and $G\subset \text{GL}\ (n+1,\C)$ a finite group acting freely off $0$ and leaving $f$ invariant.  The map $f:(\C^{n+1}/G,0)\rightarrow (\C,0)$ is a $\Q$-Gorenstein smoothing of $(V,0)\equiv (W/G,0)$, with smoothing invariants $\bar{\mu}$ and $\bar{\tau}$.  
\begin{enumerate}
\item If $G\subset \text{SL}\ (n+1,\C)$, then $(V,0)$ is Gorenstein, and $\bar{\mu}\geq \bar{\tau}$, with equality if and only if $(V,0)$ is quasi-homogeneous.
\item If $G\not \subset \text{SL}\ (n+1,\C)$, then $(V,0)$ is not Gorenstein, and $\bar{\mu}\geq \bar{\tau}-1$, with equality if and only if $(V,0)$ is quasi-homogeneous.
\end{enumerate}
\end{theorem*}
(The Milnor fibre for the quotient smoothing again has rational homology only in dimension $n$, of rank $\bar{\mu}$).  The key ingredients of the proof are the main results of \cite{wall} and \cite{gl} plus the Lefschetz fixed-point theorem.

 
A normal surface singularity in $(\C^4,0)$ is smoothable, with a smooth base space for the semi-universal deformation.  Gorenstein examples are complete intersections.  It is proved in \cite{wahls}(3.14.4) that $\alpha$ is zero.  The Conjecture in Section 6 is

\begin{Codim} Let $(V,0)$ be a normal surface singularity in $(\C^4,0)$, not a complete intersection.  Then $$\mu \geq \tau -1,$$ with equality if and only if $(V,0)$ is quasi-homogeneous.
\end{Codim}


Examples are given (\ref{okuma}) of $(V,0)\subset (\C^4,0)$ which are neither  $\Q$-Gorenstein nor quasi-homogeneous, with $\mu=\tau$ (consistent with the Conjecture).  

Our original motivation for the Conjecture (even for rational singularities) concerns those $(V,0)$ admitting a ``rational homology disk smoothing," i.e., a smoothing with Milnor number $0$.  These are especially interesting to topologists because the link of such a singularity possesses a symplectic filling with no rational homology.  A complete classification exists in case $(V,0)$ is quasi-homogeneous (\cite{SSZ}, \cite{bs}), 
 and it was conjectured in \cite{way} that these are the only examples. 
If a $\mu=0$ smoothing were $\Q$-Gorenstein (as happens in the quasi-homogeneous case \cite{wah}), one would have $\alpha=0$ and $\tau \geq 1$, whence (assuming the Rational or $\Q$-Gorenstein Conjecture) $\tau=1$ and $(V,0)$ would  be quasi-homogeneous.  Thus, there would be no other examples with $\Q$-Gorenstein smoothings.  

Note finally that there exist similar results and conjectures for a reduced curve singularity $(C,0)$.  Here, $\mu$ is defined by Buchweitz-Greuel \cite{bg}, and $\tau$ (the dimension of a smoothing component) is (via  Deligne \cite{del}) an expression in terms of curve invariants.   Work of Greuel and others (\cite{g}, \cite{gmp}) gives
\begin{enumerate} 
 \item For $(C,0)$ Gorenstein, one has $\mu \geq \tau$, with equality if and only if $(C,0)$ is quasi-homogeneous.
\item For $(C,0)$ quasi-homogeneous, one has $\mu=\tau +1-t$, where $t$ (the ``type" of $(C,0)$) is the minimal number of generators of the dualizing sheaf $\omega_C$.
\item (Conjecture)  $\mu \geq \tau +1-t$, with equality if and only if $(C,0)$ is quasi-homogeneous.
\end{enumerate}

We thank Shrawan Kumar for the proof of Lemma \ref{Lef} and Duco van Straten and Jacob Fowler for help with computer calculations.

\section{{Formulas for $\mu$ and $\tau$}}

Formulas for $\mu$ and $\tau$ for a smoothing $\pi:(\mathcal V,0)\rightarrow (\C,0)$ are given in Theorem 3.13 of \cite{wahls}, modulo three conjectures, proved there only in certain cases, but later established in general.  Specifically, the first betti number of the Milnor fibre is $0$ (\cite{gst}); any smoothing can be appropriately globalized (\cite{lo}); and the dimension of the smoothing component corresponding to $\pi$ is the length of the cokernel of $$\Theta_{\mathcal V/\C}\otimes \mathcal O_V\rightarrow \Theta_V$$
(\cite{gl}).
The invariant measuring the change in the dual of the dualizing sheaves is $$\alpha= \l(\text{Coker}\ ( \omega_{\mathcal V/\C}^*\otimes \mathcal O_V\rightarrow \omega_V^*)).$$  Writing $R=\mathcal O_{(V,0)}$ for the local ring of $V$ at $0$, one has (\cite{wahls}, Cor.A.2)  $$0\leq \alpha  \leq l(\text{Ext}^1_R(\omega,R)).$$  One deduces that $\alpha=0$ if $(V,0)$ is Gorenstein or a normal surface in $\C^4$. 

We rewrite the formulas relating $\mu, \tau,$ and $\alpha$ not in terms of $\Theta_X$ and its second exterior power $-K_X$, but rather in terms of the sheaf  $S=S_X$ of derivations logarithmic along $E$, and its second exterior power $-(K_X+E)$.  
\begin{lemma} (\cite{ok})\label{ok} Let  $(X,E)\rightarrow (V,0)$ be the minimal good resolution of a normal surface singularity, not an RDP.  Then \begin{enumerate}
\item $H^1_E(X,-(K_X+E))=0$
\item $h^1(X,-(K_X+E))$ equals the second plurigenus $$\delta_2(V)=\text{dim}\  H^0(X-E,2K_X+E)/H^0(X,2K_X+E)).$$
\end{enumerate}
\begin{proof} Corollary 1.9 of \cite{ok} asserts that $h^1(2K_X+E)=0$, whence by the long exact sequence in local cohomology $\delta_2(V)=h^1_E(2K_X+E).$  Local duality for a line bundle $L$ on $X$ yields $h^1_E(L)=h^1(K_X-L)$, from which both assertions follow.
\end{proof}
\end{lemma}
\begin{proposition}\label{t} Let $\pi:(\mathcal V,0)\rightarrow (\C,0)$ be a smoothing of a normal surface singularity $(V,0)$, not an RDP, with $\mu, \tau$, and $\alpha$ as before. Let  $(X,E)\rightarrow (V,0)$ be the minimal good resolution.  Then
\begin{enumerate}
\item $1+\mu=\alpha+13h^1(\mathcal O_X)+\chi_T(E)-(1/2)E\cdot(E+3K)-h^1(-(K_X+E))$
\medskip

\item \ \ \ $\tau=2\alpha +12h^1(\mathcal O_X)+\chi_T(E)-(1/2)E\cdot(E+3K)+h^1(S)-2h^1(-(K_X+E)).$

\end{enumerate}



\begin{proof}  For each exceptional $E_i$, denote the genus by $g_i$, the degree by $-d_i$, and the number of intersections with other curves by $t_i$.  
One has the standard short exact sequence $$0\rightarrow S\rightarrow \Theta \rightarrow \oplus N_{E_i}\rightarrow 0,$$ where $N_{E_i}$ is the normal bundle of  $E_i$.  Since $h^0(N_{E_i})=0,$ one has $h^1(N_{E_i})=g_i+d_i-1,$ so
$$h^1(\Theta)=h^1(S)+\sum (g_i+d_i-1).$$

\begin{lemma} Let $(X,E)\rightarrow(V,0)$ be the minimal good resolution of a normal surface singularity, not an RDP.  Then $$H^0(X,(-K_X)\otimes \mathcal O_E)=0.$$
\begin{proof} The group $H^1_E(-(K_X+E))$ is the direct limit of the direct system given by the injective maps  $$H^0(-(K_X+E)\otimes \mathcal O_Z(Z))\rightarrow H^0(-(K_X+E)\otimes \mathcal O_{Z+Z'}(Z+Z')),$$  where $Z$ and $Z'$ are effective exceptional divisors.  As $h^1_E(-(K_X+E))=0$, each of these spaces is $0$, in particular for $Z=E$, as desired.
\end{proof}
\end{lemma}
\begin{lemma} Let $(X,E)\rightarrow(V,0)$ be the minimal good resolution of a normal surface singularity, not an RDP.  Then $$h^1(X,-K_X)=h^1(X,-(K_X+E))+(1/2)E\cdot(E+3K).$$
\begin{proof} By the last Lemma and the standard exact sequence, it suffices to compute $h^1((-K_X)\otimes \mathcal O_E)=-\chi((-K_X)\otimes \mathcal O_E).$  Riemann-Roch for a line bundle $L$ and any exceptional divisor $Z$ states $$\chi(L\otimes \mathcal O_Z)=(-1/2)Z\cdot (Z+K)+Z\cdot L,$$ from which the result follows.
\end{proof}

\end{lemma}
Finishing the proof of Proposition \ref{t}, the last Lemma provides the formula for $1+\mu$.  The formula for $\tau$ requires checking the easily verified relation $$\chi_T(E)=\sum (g_i+d_i-1)  \ -E\cdot (E+3K)/2.$$
\end{proof}
\end{proposition}

\begin{corollary}\label{s} Notation as before, for a smoothing of a normal surface singularity (not an RDP), one has 
$$1+(\mu-\tau) +\alpha=h^1(\mathcal O_X)-h^1(S_X)+h^1(-(K_X+E)).$$
\end{corollary}

This allows one to formulate the Main Conjecture of the paper, which is motivated by (but independent of) smoothing questions.

\begin{conjecture}\label{conj} Let $(X,E)\rightarrow (V,0)$ be the minimal good resolution of a normal surface singularity, not an RDP.
\begin{enumerate}
\item If $(V,0)$ is Gorenstein, then $h^1(\mathcal O_X)-h^1(S_X)+h^1(-(K_X+E))\geq 1$, with equality if and only if $(V,0)$ is quasi-homogeneous.
\item If $(V,0)$ is not  Gorenstein, then $h^1(\mathcal O_X)-h^1(S_X)+h^1(-(K_X+E))\geq 0$, with equality if and only if $(V,0)$ is quasi-homogeneous.
\end{enumerate}
\end{conjecture}

As previously mentioned, $(1)$ has been proved in \cite{w}.

The relevant expression in the Conjecture can be written in an alternative and suggestive way.  The standard Euler characteristic of a locally free sheaf $\mathcal F$ on $X$ is
$$\chi(\mathcal F)=\text{dim}\ H^0(X-E,\mathcal F)/H^0(X,\mathcal F)\ +\text{dim}\ H^1(\mathcal F).$$  If $H^1_E(\mathcal F)=0$, one has simply $\chi(\mathcal F)=h^1(\mathcal F)$.  Now, in the case at hand, one has 
\begin{enumerate}
\item $h^1_E(\wedge^0 S_X)=h^1_E(\mathcal O_X)=0$, by Grauert-Riemenschneider
\item $h^1_E(S_X)=0$, by the main theorem of \cite{jonathan}
\item $h^1_E(\wedge^2(S_X)=0$ by \cite{ok}, as in Lemma \ref{ok}.
\end{enumerate}

\begin{proposition} \label{ch} On the minimal good resolution $X$ of a normal surface singularity, one has $$h^1(\mathcal O_X)-h^1(S_X)+h^1(-(K_X+E))=\sum_{i=0}^2 (-1)^i\chi(\wedge^i S_X).$$
\end{proposition}
It is interesting to compare with Riemann-Roch for a rank 2 vector bundle $\mathcal F$ on a smooth projective surface $Y$, which yields $$\sum_{i=0}^2 (-1)^i\chi(Y,\wedge^i \mathcal F)=c_2(\mathcal F).$$

\section{{The quasi-homogeneous case}}
The following result was asserted in \cite{wahls}(4.10.2), but not carefully proved there.

\begin{proposition}  Let $(X,E)\rightarrow (V,0)$ be the minimal good resolution of a quasi-homogeneous surface singularity, not a cyclic quotient.  Then the Euler derivation $D$ of $(V,0)$ induces a nowhere-0 section of the vector bundle $S_X$,  hence gives a short exact sequence $$0\rightarrow \mathcal O_X\rightarrow S_X\rightarrow \wedge^2(S_X)\rightarrow 0.$$
\end{proposition}
We may assume $V=\text{Spec}\ A$ is an affine variety, where $A=\oplus A_i$ is a positively graded normal domain.  Write $A=\C[z_1,...,z_s]/(g_{\alpha}(z_i))$, a quotient of a graded polynomial ring, where $\text{deg}\ z_i=m_i$.  The Euler derivation on $A_i$ is multiplication by $m_i$, and $D=\sum_{i=1}^s m_iz_i\partial/\partial z_i.$ 

Let $W\rightarrow V$ be the partial resolution obtained by blowing up the weight filtration of $A$.  $W$ has cyclic quotient singularities along the (smooth) exceptional curve $C$, which is isomorphic to Proj $A$.  Minimally resolving these cyclic quotients gives $X$ and $E$, which provides the minimal good resolution.   $D$ lifts to a section of both $\Theta_W$ and $\Theta_X$; as $H^0(S_X)=H^0(\Theta_X),$ it is also a section of $S_X$.  We must show it is nowhere $0$.  This is clear on $X-E=V-\{0\}$.

$W$ is the union of $s$ affines $W_i$, each the quotient of a smooth affine $U_i$ by a cyclic group of order $m_i$.  To define $U_1$, write $$z_1=x^{m_1}, z_2=x^{m_2}y_2,...,z_s=x^{m_s}y_s.$$  Then $U_1$ is the affine variety with coordinate ring $$\C[x,y_2,...,y_s]/(g_{\alpha}(1,y_2,...,y_s)),$$
hence is polynomial in $x$.  $D$ lifts to the derivation $x\partial/\partial x$, and is a nowhere-zero section of the bundle of derivations on $U_1$ logarithmic along $x=0.$  $W_1$ is the quotient of $U_1$ by the cyclic group of order $m_1$ generated by $$T=(1/m_1)[-1,m_2,...,m_s].$$
If the action is free at a point of $x=0$, the quotient map is a local analytic isomorphism; so at the corresponding point of $C\subset W_1$,  $D$ is still a nowhere-zero section of the corresponding bundle of logarithmic derivations.  At points of $x=0$ where there is isotropy (i.e., above the cyclic quotient singularities of $W_1$), consider the corresponding local analytic model on $U_1$.  

Changing notation slightly, for local analytic coordinates $x,y$ on $\C^2$, consider the action of $\Z/r$  of type $1/r[1,a]$, and the cyclic quotient singularity $W=\C^2/(\Z/r)$.    Let $X\rightarrow W$ be the minimal equivariant resolution of $W$, with exceptional divisor $E=\Sigma_{i=1}^l E_i$.  Denote by $C$ the Weil divisor on $W$ given by the image of $x=0$, and by $C'$ its proper transform on $Y$.  Then $C'$ intersects $E$ transversally along one end, say $E_l$.  The following Lemma will complete the proof of the Proposition.

\begin{lemma}The derivation $D=x\partial/\partial x$ on $\C^2$ induces a derivation of $W$ and lifts to a nowhere-zero section of the rank 2 vector bundle $\Omega^1_X(log(E+C'))^*$ on $X$.  

\begin{proof}  Use the familiar description (cf. Miles Reid's Warwick notes \cite{reid}, page 10) of the minimal resolution as a union $X_0 \cup X_1\cup \cdots \cup X_l$ of copies of $\C^2$.  $X_i$ has coordinates $u_i,v_i$, with the exceptional curve given by $v_0=0$ on $X_0$; $u_l=0$ on  $X_l$; and $u_iv_i=0$ on the intermediate $X_i$.   The curve $C'$ is given by $u_0=0$ on $X_0.$  On $X_i$, we have $u_i=x^{a_i}y^{b_i}, v_i=x^{c_i}y^{d_i},$ for appropriate integer exponents.  Thus on this affine, one has $D=a_iu_i\partial/\partial u_i+c_iv_i\partial/\partial v_i$.  The construction shows all the $a_i$ and $c_i$ are non-zero, except that $c_l=0$.  Thus, $D$ has the desired property on every $X_i$.
\end{proof}
\end{lemma}

\begin{theorem}\label{d} Let $(X,E)\rightarrow (V,0)$ be the minimal good resolution of a quasi-homogeneous singularity, not an RDP.
  \begin{enumerate}
 \item If $(V,0)$ is not Gorenstein, then $$h^1(\mathcal O_X)-h^1(S)+h^1(-(K_X+E))= 0.$$ 
  \item If $(V,0)$ is Gorenstein, then $$h^1(\mathcal O_X)-h^1(S)+h^1(-(K_X+E))= 1.$$ 
 \end{enumerate}
\begin{proof}  As previously noted, the Theorem is already proved for Gorenstein singularities.  It is also true for non-RDP cyclic quotients (since each individual term in $(1)$ vanishes, for instance using \ref{h1} below).
   
    For $(V,0)$ non-Gorenstein, by the preceding Proposition one must show that the (graded) coboundary map 
     $$\delta:\ H^0(X,-(K_X+E))\rightarrow H^1(X,\mathcal O_X)$$ is $0$.
  We examine the graded pieces of the map, using the notation and results of Demazure-Watanabe in \cite{wat}.  The analytic data of $(V,0)$ is given by the central curve $C=\text{Proj} \ A$, its normal divisor $D$ on $X$, and cyclic quotient singularities of type $n_i/q_i$ at points $P_i\ , i=1,...,t.$  Consider the $\Q$-divisors $$F=D-\sum (q_i/n_i)P_i$$ and $$\Xi=K+\sum (1-1/n_i)P_i.$$  Recall that for these types of $\Q$-divisors $G$, one defines an invertible sheaf 
  $$\mathcal O(G)\equiv \mathcal O(\lfloor G \rfloor))\subset k(C).$$  Pinkham's basic result is $$A=\oplus A_k=H^0(X,\mathcal O_X)=\oplus _{k=-\infty}^{\infty}H^0(C,kF)T^k.$$
  Watanabe proved \cite{wat} 
 $$ \omega_A=H^0(X-E,K_X)=\oplus _{k=-\infty}^{\infty}H^0(C,\Xi+kF)T^k.$$
It follows from the general constructions that $$H^1(X,\mathcal O_X)=\oplus _{k=0}^{\infty}H^1(C,kF)T^k.$$
$$ \omega^*_A=H^0(X-E,-K_X)=\oplus _{k=-\infty}^{\infty}H^0(C,-\Xi+kF)T^k.$$
Since $H^0(X,-(K_X+E))\subset H^0(X-E,-(K_X+E))=\omega^*_A,$ it follows that the $k$th graded piece of the coboundary map factors through $$H^0(C,-\Xi+kF)\rightarrow H^1(C,kF).$$
The second space is dual to $H^0(C,K-\lfloor kF\rfloor)$; in order for the map to be non-$0$,  both this space and $H^0(C,\lfloor -\Xi+kF\rfloor)$ must be non-0.  Therefore, the sum of the corresponding integral divisors is effective.  One checks that the sum is supported at the $P_i$'s, with coefficients 
$$\lfloor -((kq_i-1)/n_i \ +1)\rfloor-\lfloor -kq_i/n_i\rfloor.$$  But this expression equals $-1$, unless $kq_i\equiv 1\  \text{mod}\ n_i$, in which case it is $0$.  Thus, a non-zero coboundary map requires $kq_i\equiv 1\  \text{mod}\ n_i$, for all $i$.  Since the sum of the two divisors is effective and of degree $0$, and each divisor has a section, the divisors themselves would have to be trivial.  This implies in particular that $K=\lfloor kF\rfloor.$  By Watanabe's criterion, $(V,0)$ is Gorenstein.

It is not necessary to prove here the delicate fact that the coboundary map has rank $1$ in the Gorenstein case \cite{h}.
 \end{proof}
 \end{theorem}
 
\section{{Rational surface singularities}}

\begin{RatConj}  Let $(X,E)\rightarrow (V,0)$ be the minimal resolution of a rational surface singularity, not an RDP.  Then $$h^1(S_X)\leq h^1(-(K_X+E)),$$ with equality if and only if $(V,0)$ is quasi-homogeneous.
\end{RatConj}

The dimension $h^1(S_X)$ of the space of equisingular deformations is very difficult to compute (but see Example \ref{c} below).  The term  $h^1(-K_X)$ has been encountered in \cite{wahl}, (1.12), where it is proved that for $R$ the local ring of $(V,0)$,
\begin{itemize}
\item $h^1(-K_X)=l(\text{Ext}^1_R(\omega,R))$
\item $h^1(-K_X) \geq \text{mult}\ R \ - \ 3.$
\end{itemize}
As $E\cdot (E+K)=-2$, Lemma 2.3 implies that
$$h^1(-K_X)=h^1(-(K_X+E))+((-E\cdot E)\ -\ 3).$$
If $E$ is the fundamental cycle, $\text{mult}\  R =-E\cdot E$, and $h^1(-(K_X+E))$ calculates the difference in the second inequality above.  

We show $h^1(-(K_X+E))$ can be calculated from the graph of $(V,0)$.
Recall the notation and results of J. Giraud \cite{gi}.  For a line bundle $L$ on a good resolution $(X,0)\rightarrow (V,0)$ of a normal surface singularity,  define $[L]$ to be the smallest (in the usual ordering) integral divisor $D$ so that $$D\cdot E_i \leq L\cdot E_i\ ,\ \text{for all}\ i.$$  Suppose $L$ is numerically equivalent to $\sum a_iE_i$ with $a_i\in \Q$. One forms $[L]$ by first rounding up all $a_i$ to form an integral $D_1$, and inductively defining $D_{n+1}=D_n+E_j$ if one finds that $D_n\cdot E_j >L\cdot E_j.$  

\begin{theorem} (Giraud \cite{gi}) Let $L$ be a line bundle on $X$ so that $$[L] \cdot E_i\ \leq 0\ , \text{all}\ E_i.$$  Then $H^1_E(L)=0.$
\end{theorem}
 For $L$ a line bundle and $D$ an integral divisor, one has $[L-D]=[L]-D;$ in particular $$[(L-[L])]=0.$$  
 
 \begin{corollary} For any line bundle $L$ on $X$, one has $H^1_E(L-[L])=0.$
 \end{corollary}
Note $[L]$ is an effective divisor if $L\cdot E_i \leq 0$ for all $i$, or more generally if the coefficients $a_i$ as above are all $>-1$.  

\begin{proposition}  On the minimal resolution $X$ of a rational singularity, let $L$ be a line bundle for which $[L]\equiv Z$ is an effective divisor.  Then 
$$dim\  H^1(X,L)=Z\cdot(Z+K)/2\ -Z\cdot L.$$
\begin{proof}  The result is easy if $Z=0$,  so assume $Z>0.$  By duality and Corollary 5.2, one has $$0=h^1_E(\mathcal O(L-Z))=h^1(\mathcal O(K-L+Z)),$$
so that  $$0=h^1(\mathcal O(K-L+Z)\otimes \mathcal O_Z)=h^1_E(\mathcal O(K-L+Z)\otimes \mathcal O_Z).$$
Consider the short exact sequence of sheaves on $X$:
$$0\rightarrow \mathcal O(K-L)\rightarrow \mathcal O(K-L+Z) \rightarrow \mathcal O(K-L+Z)\otimes \mathcal O_Z\rightarrow 0.$$
Taking local cohomology yields $$0\rightarrow H^0_E(\mathcal O(K-L+Z)\otimes \mathcal O_Z)\rightarrow H^1_E(\mathcal O(K-L))\rightarrow H^1_E(\mathcal O(K-L+Z)) \rightarrow 0.$$
The first term (because of vanishing of the corresponding $H^1_E$) is the Euler characteristic of a line bundle on $Z$, so by Riemann-Roch its dimension is $Z\cdot(Z+K)/2\ -Z\cdot L.$  The middle term has dimension $h^1(L)$, and one must show 
$h^1_E(K-L+Z)=0$.  As usual, it suffices to prove that for any effective cycle $Y$, there is an $E_i$ in the support of $Y$ so that $$(K-L+Z+Y)\cdot E_i <0.$$  By definion $(Z-L)\cdot E_i \leq 0$ for all $i$; so it suffices to find $(K+Y)\cdot E_i<0.$  This is standard, because by rationality one always has $Y\cdot(Y+K)\leq -2.$
\end{proof}
\end{proposition}
\begin{corollary}\label{h1} Let $(X,E)\rightarrow (V,0)$ be the minimal resolution of a rational surface singularity, not an RDP.  Then
\begin{enumerate}
\item Denoting $Y=[-K_X]$, one has $$h^1(X,-K_X)=Y\cdot (Y+3K)/2.$$
\item Denoting $Z=[-(K_X+E)]=Y-E$, one has $$h^1(X,-(K_X+E))=Z\cdot(Z+3K)/2 +Z\cdot E.$$ 
\item In particular, these dimensions are topological, and depend only on the graph of $E$.
\end{enumerate}
\begin{proof}  Since $K\cdot E_i \geq 0$ for all $i$ and $K\neq 0$, one has that $-K\equiv \sum a_iE_i$, with all $a_i>0.$  Thus, both $[-K]$ and $[-(K+E)]$ are effective, so the Proposition applies.
\end{proof} 
\end{corollary}

\begin{remark}\label{r}  For a rational singularity, $[-K]\cdot E_i\leq 2-d_i\leq 0$, for all $i$, so that $[-K]$ is at least as big as the fundamental cycle $Z_0$.  The easiest way to compute $[-K]$ (and hence $[-K]-E$) is to add curves to $Z_0$ until one reaches a cycle $Y$ satisfying $Y\cdot E_i\leq 2-d_i$ for all $i$.
\end{remark}

\begin{remark} It follows from \cite{LW}, (4.11.4) that if a normal surface singularity admits a $\Q$-Gorenstein smoothing, then $h^1(-K_X)=-K\cdot K \ +h^1(\mathcal O_X).$
\end{remark}

\begin{proposition}\label{b}  Consider a rational surface singularity $(V,0)$ with a star-shaped graph (not a cyclic quotient singularity), with $(X,E)\rightarrow(V,0)$ the minimal resolution.  Then the Main Conjecture holds for $(V,0)$, i.e. $$h^1(X,S_X)\leq h^1(X,-(K_X+E)), $$
with equality iff $(V,0)$ is quasi-homogeneous.
\begin{proof}   In the local ring of $V$, consider the filtration given by the order of vanishing along the central curve in the resolution.   The associated graded ring $A$ for this filtration is normal, and is a quasi-homogeneous rational singularity with isomorphic resolution graph. The corresponding degeneration gives a deformation of Spec $A$ whose general fibre is isomorphic to $(V,0)$, and this family admits a simultaneous equisingular resolution (and in fact the reduced exceptional curve is analytically isomorphic). By the previous Corollary, $h^1(-(K+E))$ is constant in the family, and this equals $h^1(S)$ on the special fibre.  It suffices to show that if $h^1(S)$ is constant in the family, then the original singularity was already quasi-homogeneous.  Let $C$ denote the central curve, and consider on any of the resolutions the exact sequence $$H^0(S)\rightarrow H^0(S\otimes \mathcal O_C)\rightarrow H^1(S(-C))\rightarrow H^1(S)\rightarrow H^1(S\otimes \mathcal O_C)\rightarrow 0.$$
According to \cite{w}, (3.11) and (3.2), the second space has dimension $1$ and the first map is surjective if and only if the singularity is quasi-homogeneous. Since the special fibre has $$h^1(S)=h^1(S(-C))+h^1(S\otimes \mathcal O_C),$$ if this quantity remained constant in a deformation then the same would be true for each summand, and the first map in the sequence would have to remain surjective.

\end{proof}

\end{proposition}
\begin{proposition}\label{a} The Conjecture is true for rational singularities for which $h^1(S)=0$ (e.g., for the \emph{taut} singularities classified by Laufer \cite{laufert}).
\begin{proof}  One must show that a non-quasi-homogeneous $(V,0)$ with $h^1(S_X)=0$ must have $h^1(-(K_X+E))>0$.  By the preceding Proposition, it suffices to assume the graph is not star-shaped.  Since $H^1(S_X)\rightarrow H^1(\Theta_E)$ is surjective, $h^1(S_X)=0$ implies all vertices in the graph have valency at most $3$.  So, the graph contains a subgraph
   $$
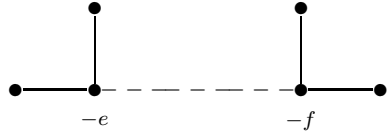
\xymatrix@R=12pt@C=24pt@M=0pt@W=0pt@H=0pt{
&&{\bullet}\lineto[dd]&&&{\bullet}\lineto[dd]\\
&&\lineto[u]&&&\lineto[u]\\
&{\bullet}\lineto[r]&\undertag
{\bullet}{-e}{0pt}\lineto[u]\dashto[r]&\dashto[r]&\dashto[r]&
\undertag{\bullet}{-f}{0pt}\lineto[r]&{\bullet}\\
&&&&
}  
$$
with all negative self-intersections at least $2$.   (Some of these can be taut.)  Let  $F$ denote the reduced curve connecting the nodes,  running from the $-e$ curve to the $-f$; it suffices to show that  $h^1(-(K_X+E)\otimes \mathcal O_F)\neq 0.$  Since $F$ is supported on a cyclic quotient singularity, the sheaf $-(K_X+E)\otimes \mathcal O_F$ depends only on the numerics of the line bundle $-(K_X+E)$ restricted to $F$.  But this line bundle dots to $0$ with the interior curves, and to $-1$ with the nodes; thus, it equals $(K_X+F)\otimes \mathcal O_F,$  whence $$h^1(-(K_X+E)\otimes \mathcal O_F)=h^1((K_X+F)\otimes \mathcal O_F).$$  Computing now on a resolution of the cyclic quotient, one has $h^1(K_X)=0$, and standard dualities and vanishing give that $$h^1((K_X+F)\otimes \mathcal O_F)=h^1(K_X+F)=h_F^1(\mathcal O (-F))=h^0(\mathcal O _F)=1.$$  \end{proof}
   \end{proposition}
   \begin{example}\label{c}Consider a rational singularity with the following graph, where as usual the unmarked bullets are $-2$ curves:
$$
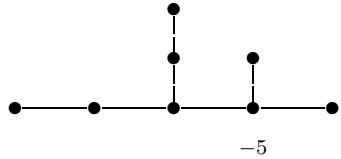
\xymatrix@R=6pt@C=24pt@M=0pt@W=0pt@H=0pt{
\\&&&&
\\
&&&\bullet\lineto[d]\\
&&&\\
&&&\bullet\lineto[u]
&\bullet\\
&&&\lineto[u]&\lineto[u]\\
&\bullet\lineto[r]&
{\bullet}\lineto[r]&\bullet\lineto[u]\lineto[r]&\undertag{\bullet}{-5}{4pt}\lineto[u]\lineto[r]
&\bullet\\
&&&&&\\
&&\\
}
$$

First, we have $h^1(-(K_X+E))=2$.  Starting from the fundamental cycle as in Remark \ref{r}, one finds $[-K_X]\equiv Y$, which has multiplicity $3$ at the left hand node, $1$ on the four outer vertices, and $2$ on the others.  Calculating with $Z\equiv Y-E$, one deduces from Corollary \ref{h1} that $h^1(-(K_X+E))=2.$  Next, consistent with the Conjecture, there is a singularity with this graph and with $h^1(S_X)=1$.  We can write it as a splice quotient \cite{nw}, i.e. the quotient of $x^3+y^3+zw^7=z^2+w^2+xy=0$ by the discriminant group $G$ (of order $60$).  Calculating the semi-universal deformation via Singular, one finds that the only deformation on which $G$ acts equivariantly is obtained by adding $tz^3w^5$ to the first equation.  By general theory, this represents exactly the only equisingular deformation of the original singularity.
\end{example}

\begin{remark} \label{jon} In a forthcoming paper, we show that the Conjecture is true for any graph satisfying $d_i\geq 2t_i-2$, for all $i$; this is a condition only at the nodes.  More precisely, if such a graph is not star-shaped, then $h^1(S_X)$ is equal to the number of ends of the graph minus $4$, which in turn is one less than $h^1(-(K_X+E)).$
\end{remark}

   %

 


\section{{ $\Q$-Gorenstein smoothings}}
Suppose $(V,0)$ is an isolated Cohen-Macaulay singularity, whose dualizing sheaf $\omega$ has order $r$.  Denoting $(\omega^{\otimes i})^{**}\equiv \omega^{[i]}$ and choosing an isomorphism $\omega^{[r]}\simeq \mathcal O$, the index one cover is the analytic spectrum of $\mathcal O \oplus \omega \oplus \omega^{[2]} \oplus \cdots \oplus \omega^{[r-1]}$.  The germ $(V,0)$ is 
said to be \emph{$\Q$-Gorenstein} if the index one cover is also Cohen-Macaulay (in which case it is Gorenstein).  (Alternatively, $(V,0)$ is the quotient of a Gorenstein singularity by a finite group acting freely off the singular point.) This is equivalent to requiring that all $\omega^{[i]}$ have maximum depth.  In particular, $\omega^*\simeq \omega^{[r-1]}$ is Cohen-Macaulay.  Recall that a rational surface singularity is $\Q$-Gorenstein.

A \emph{$\Q$-Gorenstein smoothing} $\pi:(\mathcal V,0)\rightarrow (\C,0)$ of $(V,0)$ is one which is the quotient of a smoothing of the index $1$ cover.

\begin{lemma}\label{alpha}  Let $\pi:(\mathcal V,0)\rightarrow (\C,0)$ be a $\Q$-Gorenstein smoothing of the normal surface singularity $(V,0)$.  
 Then $\alpha_{\pi}=0.$
\begin{proof}  As indicated above, the dualizing sheaf $\omega_{\mathcal V}^*\simeq \omega_{\mathcal V/\C}^*$ has depth $3$.  Tensoring with $\mathcal O_V$ gives a depth $2$ subsheaf of $\omega_V$ which is equal off the singular point, hence equal.  Therefore, $\alpha_{\pi}=0$.
\end{proof}
\end{lemma}

Suppose $G\subset \text{GL}(n+1,\C)$ is a finite group acting freely off the origin, with $n\geq 2$.  Then $(\C^{n+1}/G,0)$ is Cohen-Macaulay, and is Gorenstein if and only if $G\subset \text{SL}(n+1,\C)$ (since $G$ contains no pseudo-reflections).   Next let $f:(\C^{n+1},0)\rightarrow (\C,0)$ be a $G$-invariant function having an isolated singularity $(W,0)$ at the origin, with invariants $\mu$ and $\tau$.  The induced map $\bar{f}:(\C^{n+1}/G,0)\rightarrow (\C,0)$ is a $\Q$-Gorenstein smoothing of the quotient singularity $(W/G,0)\equiv(V,0)$.   $G$ acts freely on the Milnor fibre $M$ of $W$, so the new smoothing has Milnor fibre $M/G$, which continues to have (rational) homology only in dimensions $0$ and $n$.  Euler characteristics (and hence Milnor numbers) are related by $$1+(-1)^n\mu=|G|(1+(-1)^n\bar{\mu}).$$  Let $\bar{\tau}$ denote the dimension of the corresponding smoothing component of $V$.

\begin{theorem}\label{q}Consider the $\Q$-Gorenstein smoothing above of $(V,0)$.
 \begin{enumerate}
\item If $G\subset \text{SL}(n+1,\C)$ , then $\bar{\mu} \geq \bar{\tau}$, with equality if and only if $(V,0)$ is quasi-homogeneous.
\item  If $G\not \subset \text{SL}(n+1,\C)$, then $\bar{\mu} \geq \bar{\tau} -1$, with equality if and only if $(V,0)$ is quasi-homogeneous.
\end{enumerate}
\end{theorem}

We will prove this result via a series of Lemmas.  Denote by $J$ and $T$ the Jacobian and Tjurina algebras of $f$; the group $G$ acts on both these spaces.

\begin{lemma} For the $G$-invariant subspaces, one has $$\text{dim}\  J^G \geq \text{dim}\  T^G,$$ with equality if and only if $f$ and hence $V$ are quasi-homogeneous.
\begin{proof} The inequality is clear.  $f$ is $G$-invariant, so is in $J^G$, and is $0$ iff it is $0$ in $J$.  If $f$ is quasi-homogeneous, one may choose weights so that $G$ commutes with the $\C^*$ action.  Thus, $V$ is also quasi-homogeneous.
\end{proof}
\end{lemma}
Let $\tau'= \text{dim}\  T^G,$ and choose $G$-invariant polynomials $h_1,\cdots,h_{\tau'}$ inducing a basis of $T^G.$  The semi-universal deformation of $(W,0)$ restricts to a family $$\mathcal W=\{f+\sum_{i=1}^{\tau'}t_ih_i=0\}\subset \C^{n+1}\times \C^{\tau'}$$
$$\downarrow $$
 $$\C^{\tau'},$$ on which $G$ acts, trivially on the base and hence acting on the fibres.  This gives a smooth family $$(\mathcal W/G,0)\equiv (\mathcal V,0) \rightarrow (\C^{\tau'},0)$$ of deformations of $(V,0)$.

\begin{lemma}  The family $(\mathcal V,0) \rightarrow (\C^{\tau'},0)$ is a full smoothing component of the base space of the semi-universal deformation of $(V,0)$.  In particular, $\tau'=\bar{\tau}$.
\begin{proof} Denote the relative Kahler differentials of $f:(\C^{n+1},0)\rightarrow (\C,0)$ by $\Omega^1_f$, and its dual by $\Theta_f$; similarly one has $\Omega^1_{\bar{f}}$ and $\Theta_{\bar{f}}$.  By \cite{gl},
$\bar{\tau}$ is the length of the cokernel of the inclusion $$\Theta_{\bar{f}}\otimes \mathcal O_V \subset \Theta_V.$$  But $\Theta_V=(\Theta_W)^G$, since $W$ has an isolated normal singularity, $G$ acts freely off that point, and $\Theta_W$ is reflexive and locally free off the singular point.  Similarly, $\Theta_{\bar{f}}=(\Theta_f)^G$.  Thus, $\bar{\tau}$ is the colength of the inclusion $$\Theta_{f} ^G \otimes \mathcal O_W \subset (\Theta_W)^G.$$
Writing $R=\C\{z_1,\cdots,z_{n+1}\}$ and using the sequence $$0\rightarrow f^*(\Omega^1_{\C\{t\}}) \rightarrow \Omega^1_R\rightarrow \Omega^1_f\rightarrow 0, $$  basic homological algebra (e.g., \cite{wahls}($A.1$) with $M=\Omega^1_f$ and $t=f$) gives the exact sequence 
$$0\rightarrow \Theta_f\rightarrow \Theta_f\rightarrow \Theta_W \rightarrow \text{Ext} ^1_R(\Omega^1_f,R)\rightarrow \text{Ext} ^1_R(\Omega^1_f,R)\rightarrow T^1_W\rightarrow \text{Ext} ^2_R(\Omega^1_f,R)=0.$$
 One has a natural $G$-isomorphism  $\text{Ext} ^1_R(\Omega^1_f,R)
\cong J_f$, so taking $G$-invariants of the $G$-equivariant exact sequence $$0\rightarrow \Theta_f\otimes \mathcal O_W\rightarrow  \Theta_W\rightarrow J_f \rightarrow J_f \rightarrow T_f\rightarrow 0$$   yields that $$\bar{\tau}=\text{dim}\  T_f^G=\tau'.$$ 
\end{proof}
\end{lemma}

To compare $\mu$ to $\text{dim}\ J_f^G,$ we relate the action of $G$ on $J_f$ with the action on the cohomology of the Milnor fibre (in particular, on $H\equiv H^n(M;\C)$.)  The $n$-form $dx_1\wedge dx_2\wedge \cdots \wedge dx_{n+1}/df$ is holomorphic on the Milnor fibre, and $G$ acts on it via the determinant character.  Recall the following:

\begin{theorem} (Wall, \cite{wall}) The $G$-modules $H$ and $J_f \otimes (\text{det})$ are isomorphic, where $det$ is the determinant character.
\end{theorem}
(This result does not require that $G$ acts freely off the origin, only that it leaves $f$ invariant.)  We conclude: 
\begin{lemma} dim $J^G=\text{dim}\ H^{\nu}$, where $H^{\nu}$ is the isotypic component of $H$ corresponding to the character $\nu=(\text{det}\ )^{-1}.$
\end{lemma}

The action of $G$ on the full cohomology ring $H^*(M)$ can be deduced from the following general result.
\begin{lemma}\label{Lef} Let $X$ be a finite connected CW complex with cohomology $H^*(X)$ only in even degrees (respectively, only in odd degrees and degree $0$).  Suppose a finite group $G$ acts freely on $X$.  Then as a $G$-module, $H^*(X)$ is exactly $\chi(X)/|G|$ copies of the regular representation (respectively,the direct sum of two trivial representations and the direct sum of $-\chi(X)/|G|$ copies of the regular representation.)
\begin{proof}  Let $\phi_g:X\rightarrow X$ be the homeomorphism corresponding to $g\in G$.  Consider the Lefschetz number $$L(\phi_g)=\sum_i (-1)^i\text{Tr}\ (\phi_g^*:H^i(X)).$$$G$ acts freely, so by the Lefschetz fixed point theorem, for $g\neq e$ one has $$L(\phi_g)=0.$$Clearly, $$L(\phi_e)=L(\text{Id})=\chi(X).$$ 

When $H^*(X)$ is non-$0$ only in even degrees, one has by definition that $$L(\phi_g)=\chi(g),$$ where $\chi$ is the character corresponding to the representation of $G$ on $H^*(X)$.  Note the Euler characteristic $\chi(X)=\text{dim}\ H^*(X)$.  Thus, one has $\chi(g)=0$ for $g\neq e$ and $\chi(e)=\text{dim}\ H^*(X)$.  Standard character theory implies that the representation is a direct sum of copies of the regular representation.

When $H^*(X)$ has cohomology only in odd degrees and degree $0$, write $H^*(X)=H^0\oplus H'$, where $H'$ is a direct sum of the odd cohomology (which we assume is non-zero).  Let $\eta$ be the character corresponding to the representation of $G$ on $H'$.  Then the argument above yields $$\eta(g)=1\  \text{if}\  g\neq e$$ $$\eta(e)=\text{dim}\ H'.$$  This time, character theory implies that the representation of $G$ on $H'$ is a direct sum of a trivial representation and (dim $H'-1)/|G|$ copies of the regular representation.  As $\chi(X)=1-\text{dim}\ H'$, the claim easily follows.

\end{proof}
\end{lemma}
Back to the Milnor fibre M, we conclude that 
\begin{enumerate}
\item If $n$ is even, then $H^*(M)$ is a direct sum of $(1+\mu)/|G|$ copies of the regular representation
\item If $n$ is odd, then $H^*(M)$ is a direct sum of two trivial representations and $(\mu -1)/|G|$ copies of the regular representation.
\end{enumerate}

In particular, for $n$ even, every multiplicative character $\chi$ of $G$ occurs in $H^*(M)$ with the same multiplicity $(1+\mu)/|G|$, which as mentioned before equals $1+\bar{\mu}$.  Thus, the dimension of $H^{\chi}$ is $(1+\mu)/|G|$ unless $\chi$ is the trivial character, in which case the dimension is $1$ less.  But the determinant character of $G$ is trivial if and only if $G\subset \text{SL}\ (n+1,\C)$.  One now has all the ingredients to verify Theorem 5.2 in this case.

When $n$ is odd, Lemma \ref{Lef}  implies that as a $G$-module, $H$ consists of one copy of the trivial representation plus $(\mu -1)/|G|=\bar{\mu}-1$ copies of the regular representation.  Thus, one knows the dimension of each isotypic component, and the same argument as in the even case completes the proof of the Theorem.

\section{{normal surfaces in $\C^4$}}

Let $(V,0)\subset (\C^4,0)$ be a germ of a normal surface singularity.  If $(V,0)$ is Gorenstein, it is a complete intersection; otherwise, it is maximal-minor determinantal.  The base space of the semi-universal deformation is smooth, as $T^2_{(V,0)}=0;$ computer programs can be used to calculate it.  It is unknown whether the Milnor fibre is always simply-connected; this is true for the rational triple points, since one has a unique Milnor fibre, diffeomorphic to the minimal resolution.  The Milnor number is difficult to compute.  Since $\alpha=0$ in these cases, one wants to prove

\begin{Codim} Let $(V,0)\subset (\C^4,0)$ be a non-Gorenstein singularity.  Then $$\mu \geq \tau -1,$$ with equality if and only if $(V,0)$ is quasi-homogeneous.
\end{Codim}

\begin{example}\label{okuma}\cite{okuma},(6.3).  There is a $(V,0)\subset (\C^4,0)$ which is not $\Q$-Gorenstein and not quasi-homogeneous (though having a star-shaped graph), with $\mu = \tau$.

Consider the singularity $(V,0)$ defined by the vanishing of the $2\times 2$ minors of the matrix
 \[\left( \begin{array}{ccc} x & y & z \\ y-3w^2  & z+w^{3m} & x^2+6wy-2w^3 \\
    \end{array} \right), \ m\geq 1. \] 
$(V,0)$ has multiplicity $3$ with $p_g=m$, and the same (integral homology sphere) resolution graph as the hypersurface singularity $x^2+y^3+z^{6m+7}=0$ (which has $p_g=m+1$):
$$
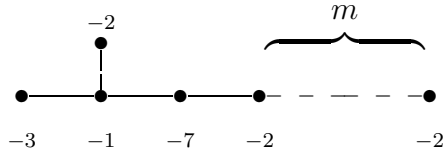
\xymatrix@R=6pt@C=24pt@M=0pt@W=0pt@H=0pt{
\\&&&&&m\\
&&\overtag{\bt}{-2}{8pt}
&&&{\hbox to 0pt{\hss$\overbrace{\hbox to 60pt{}}$\hss}}&\\
&&\lineto[u]\\
&\undertag{\bt}{-3}{6pt}\lineto[r]&\undertag{\bt}{-1}
{6pt}\lineto[u]\lineto[r]&\undertag{\bt}{-7}{6pt}\lineto[r]&\undertag{\bt}{-2}{6pt}\dashto[r]&\dashto[r]
&\undertag{\bt}{-2}{6pt}\\
&&&&\\
&&&&\\
}
$$

\vspace{.2in}

Corollary $6.6$ of \cite{okuma} yields that $\mu=12m+1.$  Duco van Straten has used Singular to prove for $m=1$, and to indicate the likelihood for general $m$, that $\tau=12m+1$.  The semi-universal deformation should be obtained by perturbing the matrix by adding
 \[\left( \begin{array}{ccc} 0 & 0 & 0 \\ f  & g & h \\
    \end{array} \right),  \] 
    where $f=\sum_{i=0}^{3m-1}  a_i w^i $ ,
$g=\sum_{i=0}^{3m-1} b_i w^i$,
$h=\sum_{i=0}^{3m}c_i w^i +x \sum_{i=0}^{3m-1} d_i w^i$.  Thus, the $a_i, b_i, c_i, d_i$ are parameters for the semi-universal deformation.
\end{example}

\bigskip

\end{document}